\documentclass[12pt]{amsart}

\usepackage{amsmath}
\usepackage{amsthm}
\usepackage{stmaryrd}
\usepackage{amsfonts,mathrsfs,amssymb}
\usepackage{fullpage}

\title[Division Theorems]{Division Theorems and Twisted Complexes}
\author{Dror Varolin}
\thanks{Partially supported by NSF grant DMS-0400909}

\address{Department of Mathematics \newline
\indent Stony Brook University \newline
\indent Stony Brook, NY 11794-3651}

\email{dror@math.sunysb.edu}

\newcommand{\noi}{\noindent}

\newcommand{\co}{{\mathcal O}}

\newcommand{\sC}{{\mathscr C}}
\newcommand{\sd}{{\mathscr D}}

\newcommand{\sF}{{\mathscr F}}

\newcommand{\sh}{{\mathscr H}}
\newcommand{\si}{{\mathscr I}}
\newcommand{\sj}{{\mathscr J}}
\newcommand{\sk}{{\mathscr K}}

\newcommand{\fD}{{\mathfrak D}}

\newcommand{\vp}{\varphi} 
\newcommand{\ve}{\varepsilon}

\newcommand{\C}{{\mathbb C}}

\newcommand{\R}{{\mathbb R}}

\newcommand{\di}{\partial}
\newcommand{\dbar}{\bar \partial}

\newcommand{\re}{{\rm Re\ }}

\newcommand{\ii}{\sqrt{-1}}

\newcommand{\tensor}{\otimes}

\begin{document}
\maketitle

\theoremstyle{plain}
\newtheorem{m-thm}{\sc Theorem}

\newtheorem{thm}{\sc Theorem}[section]
\newtheorem{lem}[thm]{\sc Lemma}
\newtheorem{prop}[thm]{\sc Proposition}
\newtheorem{cor}[thm]{\sc Corollary}
\newtheorem{defn}[thm]{\sc Definition}

\theoremstyle{definition}
\newtheorem{conj}[thm]{\sc Conjecture}
\newtheorem{qn}[thm]{\sc Question}
\newtheorem{pr}[thm]{\sc Problem}

\theoremstyle{remark}
\newtheorem*{rmk}{\sc Remark}
\newtheorem*{ack}{\sc Acknowledgment}
\newtheorem*{ex}{\sc Example}

\section{Introduction}

A classical problem both in commutative algebra and in several complex variables is the ideal membership, or division problem.  In the setting of commutative algebra, the decisive result is  Hilbert's Nullstellensatz.  On the other hand, in several complex variables the most basic division problem for bounded holomorphic functions is open:  given a collection of bounded holomorphic functions $g^1,...,g^p$ on the unit ball $B \subset \C^n$ ($n \ge 2$) such that $\sum |g^i|^2 = 1$, are there bounded holomorphic functions $h_1,...,h_p$ such that $\sum h_i g^i = 1$?  The case $n=1$ is the famous Corona Theorem of Carleson.  In higher dimensions, the closest result thus far is an $L^2$ version, also known as the celebrated Division Theorem of Skoda.

In this paper, we use the method of the twisted Bochner-Kodaira Identity together with Skoda's Basic Estimate to establish a generalization of Skoda's Division Theorem.

Skoda's Theorem has many applications.  Some recent applications in algebraic geometry appear in the work of Ein and Lazarsfeld \cite{el}, where Skoda's Theorem plays a key role in giving an effective version of the Nullstellensatz.  Siu has used Skoda's Theorem to prove results about effective global generation of multiplier ideals, which was a key tool in his approach to establishing the deformation invariance of plurigenera \cite{siu-pluri}.  Siu has also used Skoda's Theorem in his approach to the problem of finite generation of the canonical ring \cite{siu-mult}.

At the same time, several authors have been establishing results that show the fundamental role of the Ohsawa-Takegoshi extension theorem and its variants in the areas of analytic methods in algebraic geometry and of several complex variables.  The applications of the Ohsawa-Takegoshi Theorem are too numerous to mention in this introduction.

The results of Skoda and Ohsawa-Takegoshi are similar in nature and proof.   Both results use the Bochner-Kodaira Identity.  However, in the Ohsawa-Takegoshi technique, the Bochner-Kodaira Identity is "twisted".

In Skoda's Theorem,  a functional analysis argument is used that is similar to the well-known Lax-Milgram Lemma.  We recall, perhaps with slight modification, Skoda's functional analysis in Section \ref{background-section}.  As usual, the functional analysis requires us to establish an {\it a priori} estimate.  This estimate is obtained from the Bochner-Kodaira Identity together with a non-trivial and very sharp inequality due to Skoda, referred to in this paper as Skoda's Inequality.  The resulting {\it a priori} estimate is referred to here as Skoda's Basic Estimate.

The main idea of this paper is to introduce twisting into Skoda's Basic estimate.  More precisely, we twist the Bochner-Kodaira Identity before applying Skoda's Inequality.  The result is a series of divison theorems whose estimates are different from the original result of Skoda.

Because many recent applications of $L^2$ theorems have been to complex and algebraic geometry, we will state our results in the more general language of singular metrics and sections of holomorphic line bundles on so-called {\it essentially Stein} manifolds:  a K\" ahler manifold $X$ is said to be essentially Stein if there exists a complex subvariety $V \subset X$ such that $X-V$ is a Stein manifold.  For example, $X$ could be a Stein manifold, in which case we can take $V = \emptyset$, or $X$ could be a smooth projective variety, in which case $V$ could be the intersection of $X$ with a hyperplane in some projective space in which $X$ is by hypothesis embedded.  A third interesting class of examples is a holomorphic family of algebraic manifolds fibered over the unit disk or over a more general Stein manifold.

We fix on $X$ two holomorphic line bundles $E \to X$ and $F\to X$, with singular metrics $e^{-\vp _E}$ and $e^{-\vp _F}$ respectively, and suppose given a collection of sections 
\[
g^1,..., g^p \in H^0(X,E)
\]
where $p \ge 1$ is some integer.

Taking cue from Skoda \cite{s-72}, we seek to determine which sections $f \in H^0(X, F+K_X)$ can be divided by $g=(g^1,..., g^p)$, in the sense that there exist sections 
\[
h_1,..., h_p \in H^0(X, F-E +K_X)
\]
satisfying the equality
\[
f = \sum _{i=1} ^p h_i g^i.
\]
Moreover, if $f$ satisfies some sort of $L^2$ estimate, what can we say about estimates for $h_1,..., h_p$?  We shall refer to this question as the {\it division problem}.

\medskip

To state our main result, it is useful to introduce the following definition.  
\begin{defn}\label{skoda-triple-defn}
For a triple $(\phi, F, q)$ where $\phi : [1,\infty) \to \R $, $F:[1,\infty ) \to [0,\infty) $ are $\sC ^2$ functions and $q >0$ is an integer, define the auxiliary function 
\[
\tau(x) = x+F(x).
\]
We call $(\phi , F, p)$ a \emph{Skoda Triple} if
\begin{equation}\label{skoda}
\tau (x) \phi '(x) + 1+F'(x) \ge 0\quad \text{and} \quad (\phi ''(x) + F''(x))\le 0.
\end{equation}
Given a Skoda Triple $(\phi , F,q)$, we define
\[
B(x)= 1+ \frac{(\tau(x)\phi '(x) +1+F'(x))}{q\tau (x)}\quad \text{ and}\quad A(x) := \left \{ 
\begin{array}{c@{\quad}c}
\frac{(1+F'(x))^2}{-(F''(x)+\phi ''(x))} & 1+ F' \not \equiv 0 \\
0 & 1+ F' \equiv 0
\end{array}
\right . .
\]
\end{defn}
\begin{rmk}
Note that 
\[
\frac{B}{\tau (B-1)} = \frac{q\tau + \tau \phi ' +1 +F'}{\tau (\tau \phi '+1+F')} = 1+ \frac{q}{\tau (\tau \phi '+1+F')}
\]
\end{rmk}

Our main result is the following.
\begin{m-thm}\label{main}
Let $X$ be an essentially Stein manifold of complex dimension $n$, $F,E\to X$ holomorphic line bundles with singular Hermitian metrics $e^{-\psi}$ and $e^{-\eta}$ respectively, and $g^1,...,g^p \in H^0(X,E)$.  Let 
\[
q = \min (p-1,n) \quad \text{and} \quad \xi = 1- \log (|g|^2e^{-\eta}).
\]
Let $(\phi , F, q)$ be a Skoda triple.  Set $\tau = \tau (\xi)$, $A=A(\xi)$ and $B=B(\xi)$, with $\tau (x)$, $A(x)$ and $B(x)$ as in Definition \ref{skoda-triple-defn}.  Assume that 
\[
\ii \di \dbar \psi \ge Bq \ii \di \dbar \eta.
\]
If $\tau$ non-constant, or if $\phi$ cannot be defined on $\R$ and still satisfy \eqref{skoda} of Definition \ref{skoda-triple-defn}, assume further that 
\[
|g|^2e^{-\eta} := \sum _{j=1} ^p |g^j|^2e^{-\eta} < 1.
\]
Then for every section $f \in H^0(X,F+K_X)$ such that 
\[
\int _X \frac{B}{\tau (B-1)}\frac{|f|^2e^{\phi (\xi)}e^{-\psi}}{(|g|^2e^{-\eta})^{q+1}} < +\infty
\]
there exist sections  $h_1,...,h_p \in H^0(X,F-E+K_X)$ such that 
\[
\sum _{j=1} ^p h_j g^j = f
\]
and 
\[
\int _X \frac{|h|^2e^{\phi (\xi)}e^{-(\psi-\eta)}}{(\tau +A)(|g|^2e^{-\eta})^q} \le \int _X \frac{B}{\tau (B-1)}\frac{|f|^2e^{\phi (\xi)}e^{-\psi}}{(|g|^2e^{-\eta})^{q+1}}.
\]
\end{m-thm}

\medskip

Theorem \ref{main} implies a large number of division theorems.  However, as stated, Theorem \ref{main} does not give a solution to the division problem unless we input a Skoda triple (always with $q=\min(n,p-1)$).  By choosing different Skoda triples $(\phi , F, q)$ (with the same $q= \min (n,p-1)$ from Theorem \ref{main}), one can obtain numerous division theorems as corollaries of Theorem \ref{main}.  In the next section we establish some of these corollaries, and hope the reader is convinced that such corollaries are easy to come by, or equivalently, that Skoda triples are easy to find.

\tableofcontents

\section{Corollaries of Theorem \ref{main}}

\begin{ex}
Fix $\alpha > 0$.  Let $q = \min (n,p-1)$, $F(x) = 1-x$ and $\phi (x) = (\alpha -1)q(x-1)$.  
The $\tau =1$ so by definition $A=0$, and $B = \alpha$.  Thus we obtain the following geometric reformulation of the famous theorem of Skoda.
\begin{thm}\label{skoda-thm-1}
Let $X$ be an essentially Stein manifold of complex dimension $n$, $F,E\to X$ holomorphic line bundles with singular Hermitian metrics $e^{-\psi}$ and $e^{-\eta}$ respectively, and $g^1,...,g^p \in H^0(X,E)$.  Assume that 
\[
\ii \di \dbar \psi \ge \alpha q \ii \di \dbar \eta.
\]
Then for any $f \in H^0(X,K_X+F)$ such that 
\[
\int _{X} \frac{|f|^2 e^{-\psi}}{(|g|^2e^{-\eta})^{\alpha q + 1}} < +\infty 
\]
there are $p$ sections $h_1,...,h_p \in H^0 (X,K_X+F-E)$ such that
\[
h_k g^k = f \quad \text{and} \quad \int _{X} \frac{|h|^2 e^{-(\psi-\eta)}}{(|g|^2e^{-\eta})^{\alpha q }} \le \frac{\alpha }{\alpha -1}\int _{X} \frac{|f|^2 e^{-\psi}}{(|g|^2e^{-\eta})^{\alpha q + 1}}.  
\]
\end{thm}

\begin{rmk}
The proof of Theorem \ref{skoda-thm-1} that we give essentially reduces to Skoda's original proof since, as we shall see, when $\tau$ is constant the twisting of the $\dbar$-complex becomes trivial. 
\end{rmk}
\end{ex}

\begin{ex}
In his paper \cite{siu-mult}, Siu derived from Skoda's Theorem the following result in the case of algebraic manifolds.  We give a slightly different proof of Siu's result, in the more general setting of almost Stein manifolds.

\begin{thm}\label{siu-skoda}
Let $X$ be an almost Stein manifold of complex dimension $n$, $L\to X$ a holomorphic line bundle, and $H\to X$ a holomorphic line bundle with non-negatively curved singular Hermitian metric $e^{-\vp}$.  Let $k \ge 1$ be an integer and fix sections $G^1,...,G^p \in H^0(X,L)$.  Define the multiplier ideals 
\[
\sj _{k+1} = \si (e^{-\vp}|G|^{-2(n+k+1)}) \quad \text{and} \quad \sj _k = \si (e^{-\vp}|G|^{-2(n+k)}).
\]
Then 
\[
H^0(X,((n+k+1)L+H+K_X) \tensor \sj _{k+1}) = \bigoplus _{j=1} ^p G_j H^0(X,((n+k)L+H+K_X) \tensor \sj_k).
\]
\end{thm}

\begin{proof}
By taking $G^{p}= ... G^n = 0$, we may assume that $q = n$.  Take $\alpha = (n+k)/n$, so that $\alpha q= n+k$.  We are going to use Theorem \ref{skoda-thm-1} with $F=(n+k+1)L+H$, $E=L$ and $g^i = G^i$.  Fix a metric $e^{-\eta}$ for $L$ having non-negative curvature current,  (for example, one could take $\eta = \log |G^1|^2$)  and let $\psi = \vp +(n+k+1)\eta$.  Then 
\[
\ii \di \dbar \psi - \alpha q \ii \di \dbar \eta = \ii \di \dbar \vp + \ii \di \dbar \eta \ge 0.
\]
Suppose $f \in H^0(X,\co _X ((n+k+1)L+H+K_X) \tensor \sj_{k+1})$.  By Theorem \ref{skoda-thm-1} there exist $h_1,...,h_p \in H^0(X,\co _X ((n+k)L+H+K_X))$ such that $h_kG^k = f$.  Moreover, 
\begin{eqnarray*}
\int _X \frac{|h|^2e^{-\vp}}{|G|^{2(n+k)}}&=& \int _X \frac{|h|^2e^{-(\psi-\eta)}}{(|g|^2e^{-\eta})^{n\alpha}} \\
&\le& \frac{\alpha}{\alpha-1} \int _X \frac{|f|^2e^{-\psi}}{(|g|^2e^{-\eta})^{n\alpha +1}}\\
& = &\frac{n+k}{n+k-1} \int _X \frac{|f|^2e^{-\vp}}{|G|^{2(n+k+1)}} < +\infty.
\end{eqnarray*}
Thus $h_i \in \sj_k$ locally (and much more).  The proof is complete.
\end{proof}
\end{ex}

The rest of the results we present here were motivated in part by our desire to say something about the case $k=0$ in Theorem \ref{siu-skoda}.  Of course, Theorem \ref{siu-skoda} as stated is not true in general if $k=0$.  But at the end of this section, we will state and prove a result similar to Theorem \ref{siu-skoda} (namely Theorem \ref{siu-skoda-k=0}) that does address the question of when sections of $(n+1)L+H+K_X$ are expressible in terms of sections of $nL+H+K_X$.

\medskip

The next example again does not make use of twisting.

\begin{ex}
For any $\ve > 0$ and any positive integer $q$ the triple $(\phi , F,q)$ where $\phi (x) = \ve \log x$ and $F(x) =1-x$ is a Skoda Triple.  Indeed, $\phi '(x) = \ve/x$ and $\phi ''(x) = - \ve /x^2 < 0$, while $1+F'(x) = 0$ and $F''(x) \equiv 0 \le 0$, and Thus (1) holds.  Furthermore, 
\[
A = 0, \quad B = 1+ \frac{\ve}{q\xi} \quad \text{and}\quad \frac{B}{\tau (B-1)} = 1+\frac{q\xi}{\ve},
\]
and thus 
\[
Bq \le q+\ve, \quad \text{while} \quad \frac{B}{\tau(B-1)} \le \frac{q+\ve}{\ve}\xi .
\]
Applying Theorem \ref{main} to the Skoda Triple $(\ve \log x, 0,q)$ where $q= \min (n,p-1)$ gives us the following result.
\begin{thm}\label{untwisted-thm-1}
Let the notation be as in Theorem \ref{main}.  Fix $\ve > 0$ and assume that $|g|^2e^{-\eta}<1$ on $X$.  Set $q = \min (n,p-1)$ and $\xi = 1- \log (|g|^2e^{-\eta})$.   Suppose that
\[
\ii \di \dbar \psi \ge(q+\ve) \ii\di \dbar \eta.
\]
Then for every $f \in H^0(X, F+K_X)$ such that
\[
\int _{X} \frac{|f|^2 e^{-\psi}\xi ^{1+\ve}}{(|g|^2e^{-\eta})^{q+1}} < +\infty,
\]
there exist sections $h_1,...,h_p \in (H^0(X, F-E+K_X)$ such that 
\[
\sum _{j=1} ^p h_j g^j = f \quad \text{and} \quad \int _{X} \frac{|h|^2 e^{-(\psi-\eta)}\xi^{\ve}} {(|g|^2e^{-\eta})^q} \le \frac{q+\ve}{\ve} \int _{X} \frac{|f|^2 e^{-\psi}\xi ^{1+\ve}}{(|g|^2e^{-\eta})^{q+1}}.
\]
\end{thm}
\end{ex}

By allowing twisting, we can improve the estimates of Theorem \ref{untwisted-thm-1}, at the cost of a little more curvature from $e^{-\psi}$.

\begin{ex}
For any $\ve > 0$ and any positive integer $q$ the triple $(\phi , F,q)$ where $\phi (x) = \ve \log x$ and $F(x) \equiv 0$ is a Skoda Triple.  Indeed, $x\phi '(x) = \ve$ and $\phi ''(x) = - \ve /x^2 < 0$, while $1+F'(x) = 1 > 0$ and $F''(x) \equiv 0 \le 0$, and Thus (S1) holds.  Furthermore, $A(x) = \frac{\tau (x) }{\ve}$, and thus (S2) holds.  Moreover
\[
A(x) = \frac{\tau (x)}{\ve}, \quad B(x) = \frac{q\tau (x) +1+\ve}{q\tau (x)} \quad \text{and}\quad \frac{\tau(B-1)}{B} = \frac{(1+\ve)}{q + \frac{1+\ve}{\tau}},
\]
and thus 
\[
Bq \le q+1+\ve, \quad \text{while} \quad \frac{B}{\tau(B-1)} \le \frac{q+1+\ve}{1+\ve}.
\]
Applying Theorem \ref{main} to the Skoda Triple $(\ve \log x, 0,q)$ where $q= \min (n,p-1)$ gives us the following result.
\begin{thm}\label{twisted-thm-1}
Let the notation be as in Theorem \ref{main}.  Fix $\ve > 0$ and assume that $|g|^2e^{-\eta}<1$ on $X$.  Set $q = \min (n,p-1)$ and $\xi = 1- \log (|g|^2e^{-\eta})$.   Suppose that
\[
\ii \di \dbar \psi \ge(q+1+\ve) \ii\di \dbar \eta.
\]
Then for every $f \in H^0(X, F+K_X)$ such that
\[
\int _{X} \frac{|f|^2 e^{-\psi}\xi ^{\ve}}{(|g|^2e^{-\eta})^{q+1}} < +\infty,
\]
there exist sections $h_1,...,h_p \in (H^0(X, F-E+K_X)$ such that 
\[
\sum _{j=1} ^p h_j g^j = f \quad \text{and} \quad \int _{X} \frac{|h|^2 e^{-(\psi-\eta)}\xi^{\ve}} {(|g|^2e^{-\eta})^q \xi} \le \frac{q+1+\ve}{\ve} \int _{X} \frac{|f|^2 e^{-\psi}\xi ^{\ve}}{(|g|^2e^{-\eta})^{q+1}}.
\]
\end{thm}
\end{ex}

\begin{rmk}
Note that if, for example, $X$ is a bounded pseudoconvex domain in $\C ^n$ and we take $E=\co$ and $\eta \equiv 0$, then Theorem \ref{twisted-thm-1} is a strict improvement over Theorem \ref{untwisted-thm-1}.  Thus twisting can sometimes get us stronger results.
\end{rmk}

Though less general, Theorems \ref{skoda-thm-1}, \ref{untwisted-thm-1} and \ref{twisted-thm-1} are aesthetically more pleasing than Theorem \ref{main}, because the integrands in the conclusions of the former are more natural.  Another method for obtaining such natural integrands is through the use of the notion of {\it denominators} introduced in \cite{mv} by McNeal and the author.

\begin{defn}\cite{mv}\label{denom-defn}
Let $\sd$ denote the class of functions $R:[1,\infty) \to [1,\infty)$ with the following properties.

\begin{enumerate}

\item[(D1)]  Each $R \in \sd$ is continuous and increasing.

\item[(D2)]  For each $R \in \sd$ the improper integral 
\[
C(R) := \int_1^\infty \frac{dt}{R(t)}
\]
is finite.
\end{enumerate}
For $\delta >0$, set
\[
G_\delta(x)=\frac 1{1+\delta}\left( 1+\frac\delta{C(R)}\int_1^x \frac{dt}{R(t)}\right),
\]
and note that this function takes values in $(0,1]$. Let
\[
F_{\delta} (x) := \int _1 ^x \frac{1-G_{\delta}(y)}{G_{\delta} (y)} dy.
\]
\begin{enumerate}
\item[(D3)] For each $R \in \sd$ there exists a constant $\delta >0$ such that
\[
K_{\delta}(R) := \sup _{x \ge 1} \frac{x+F_{\delta}(x)}{R(x)}
\]
is finite.
\end{enumerate}
A function $R \in \sd$ is called a {\it denominator}.
\end{defn}

The following key lemma about denominators was proved in \cite{mv}.

\begin{lem}\label{ode-lem} If $R\in\fD$ and $\int _1 ^{\infty} \frac{dt}{R(t)} = 1$, then the function $F = F_{\delta}$ given in Definition \ref{denom-defn} satisfies
\begin{subequations}\label{h-conditions}
\begin{align}
& x+F(x) \ge 1,\\
& 1+ F'(x) \ge 1,\ \text{ and} \\
& F''(x) < 0.
\end{align}
\end{subequations}
Moreover, $F$ satisfies the ODE
\begin{equation}\label{ode}
F''(x) +\frac {\delta}{(1+\delta) R(x)}\left(1+F'(x)\right)^2 =0,\qquad x\geq 1,
\end{equation}
where $\delta$ is a positive number guaranteed by Condition {\rm (D3)} of Definition \ref{denom-defn}. 
\end{lem}
Using Lemma \ref{ode-lem} we can prove the following theorem.
\begin{thm}\label{main-1}
Let $X$ be an essentially Stein manifold of complex dimension $n$, $F,E\to X$ holomorphic line bundles with singular Hermitian metrics $e^{-\psi}$ and $e^{-\eta}$ respectively, and $g^1,...,g^p \in H^0(X,E)$.  Let 
\[
q = \min (p-1,n) \quad \text{and} \quad \xi = 1- \log (|g|^2e^{-\eta}).
\]
Let $R \in \sd$ with constant $\delta > 0$ and function $F=F_{\delta}$ determined by Definition \ref{denom-defn}, and set $B= 1+ \frac{1+F'(\xi)}{q(\xi +F(\xi))}$.  Assume that 
\[
\ii \di \dbar \psi \ge qB \ii \di \dbar \eta \quad \text{and}\quad |g|^2e^{-\eta} := \sum _{j=1} ^p |g^j|^2e^{-\eta} < 1.
\]
Then for every section $f \in H^0(X,F+K_X)$ such that 
\[
\int _X \frac{|f|^2e^{-\psi}}{(|g|^2e^{-\eta})^{q+1}} < +\infty
\]
there exist sections  $h_1,...,h_p \in H^0(X,F-E+K_X)$ such that 
\[
\sum _{j=1} ^p h_j g^j = f
\]
and 
\[
\int _X \frac{|h|^2e^{-(\psi-\eta)}}{(|g|^2e^{-\eta})^qR(\xi)} \le (1+q) \left (\frac{(1+\delta)}{\delta} C(R) + K_{\delta}(R)\right ) \int _X \frac{|f|^2e^{-\psi}}{(|g|^2e^{-\eta})^{q+1}}.
\]
\end{thm}
\begin{rmk}
For the reader that does not like the appearance of the function $B$ in the statement of Theorem \ref{main-1}, we note that $qB$ is always bounded above by $q+2+\delta$.  Indeed, note that $qB =q+ \frac{1+F'}{\tau}$.  Now, in the notation of Definition \ref{denom-defn},
\[
\tau  - (1+F'(x)) = \int _1 ^x \frac{dy}{G_{\delta}(y)} - \frac{1}{G_{\delta}(x)} \ge -(1+\delta)
\]
since, from the definition of $G_{\delta}$, $\frac{-1}{G_{\delta}} \ge -(1+\delta)$ while $\frac{1}{G_{\delta}} \ge 0$.  Thus 
\begin{equation}\label{qb-bd}
qB \le q+1+\frac{1+\delta}{\tau} \le q+2+\delta.
\end{equation}
The reason we did not hypothesize that $\ii \di \dbar \psi \ge (q+2+\delta) \ii \di \dbar \eta$ is that often can do better than the bound \eqref{qb-bd}.  It is often easy to estimate $\frac{1+F'}{\tau}$ in specific examples.
\end{rmk}

\begin{proof}[Proof of Theorem \ref{main-1}]
Let $F$ be the function associated to the denominator $R$ via Lemma \ref{ode-lem}.  Observe that in view of Lemma \ref{ode-lem}, $(0,F,q)$ is a Skoda triple.  Let $\tau$, $A$ and $B$ be the functions associated to $(0,F,q)$ in Definition \ref{skoda-triple-defn}.  We claim that 
\begin{equation}\label{ode-cor}
\frac{\tau +A}{R} \le \frac{(1+\delta)}{\delta} C(R) + K_{\delta}(R).
\end{equation}
Indeed, $\tau /R \le K_{\delta}(R)$ by property (D3), while $A/R \le \frac{(1+\delta)}{\delta} C(R)$ by the definition of $A$ and the ODE of Lemma \ref{ode-lem}.  (Note that if $R \in \sd$, then $C(R) R \in \sd$ and $\int _1 ^{\infty}(C(R)R(t))^{-1} dt = 1$.)

Finally, observe that
\begin{equation}\label{b-ineq}
\frac{B}{\tau (B-1)} =  1+\frac{q}{\tau(1+F')} \le 1+q.
\end{equation}
Thus 
\begin{eqnarray*}
&& \int _X \frac{|h|^2e^{-(\psi-\eta)}}{(|g|^2e^{-\eta})^qR(\xi)} \\
&\le& \int _X \frac{\tau+A}{R(\xi)} \frac{|h|^2e^{\phi (\xi)}e^{-(\psi-\eta)}}{(\tau +A)(|g|^2e^{-\eta})^q}\\
&\le& \left (\frac{(1+\delta)}{\delta} C(R) + K_{\delta}(R)\right )\int _X \frac{B}{\tau (B-1)}\frac{|f|^2e^{\phi (\xi)}e^{-\psi}}{(|g|^2e^{-\eta})^{q+1}}\\
&\le &(1+q) \left (\frac{(1+\delta)}{\delta} C(R) + K_{\delta}(R)\right ) \int _X \frac{|f|^2e^{-\psi}}{(|g|^2e^{-\eta})^{q+1}}.
\end{eqnarray*}
In going from the second to the third line, we used Theorem \ref{main} and The inequality \eqref{ode-cor}, and in going from the third line to the last we used \eqref{b-ineq}.  The proof is complete.
\end{proof}

The following is a table of denominators and their corresponding constants \cite{mv}.
\begin{equation}\label{table}
\begin{array}{|cc|c|}
\hline
&& \\
&R(x) = & \frac{(1+\delta)}{\delta} C(R) + K_{\delta}(R) \le \\
&&\\
\hline
&& \\
(i)& e^{s(x-1)} & \frac{4}{s} \\
&& \\
\hline
&& \\
(ii)&x^2 & \frac{3+2\sqrt{2}}{4} \\
&& \\
\hline
&& \\
(iii)&x^{1+s} & \frac{4}{s} \\
&& \\
\hline
&& \\
(iv)&R_N(x) & \frac{4}{s} \\
&& \\
\hline
\end{array}
\end{equation}
In entry $(iv)$, 
\[
R_N(x) =  x \left ( \prod _{j=1}^{N-2} L_j (x) \right ) (L_{N-1}(x))^{1+s}, 
\]
where 
\[
E_j = \exp ^{(j)}(1) \quad \text{and}\quad L_j(x) =\log ^{(j)}(E_jx).
\]

\begin{rmk}
To define denominators yielding Skoda triples $(\phi , F, q)$ with both $\phi$ and $F$ non-trivial seems a little more complicated.  The corresponding ODE that determines the associated function $\tau$ does have solutions, but since this (first order) ODE is not autonomous, it is harder to get explicit properties of the associated function $\tau$.
\end{rmk}

We can now state and prove our Siu-type division theorem for the case $k=0$, under an additional assumption on the line bundle $L$.

\begin{thm}\label{siu-skoda-k=0}
Let $X$ be an almost Stein manifold of complex dimension $n$, $L\to X$ a holomorphic line bundle, $H\to X$ a holomorphic line bundle with non-negatively curved singular Hermitian metric $e^{-\vp}$.  Fix sections $G^1,...,G^p \in H^0(X,L)$ and a singular Hermitian metric $e^{-\eta}$ for $L$ having non-negative curvature, and such that 
\[
|G|^2e^{-\eta} < 1 \quad \text{on }X.
\]
Fix a denominator $R \in \sd$ such that the associated function $B$ satisfies $nB \le n+1$.  (For example, this is the case for the denominator $(i)$ of Table \eqref{table}.)  Define the multiplier ideals 
\[
\si _{1} = \si \left (\frac{e^{-\vp}}{|G|^{2(n+1)}} \right ) \quad \text{and} \quad \si _0 = \si \left (\frac{e^{-\vp}}{|G|^{2n}R(1-\log |G|^2 +\eta)}\right ).
\]
Then 
\[
H^0(X,((n+1)L+H+K_X) \tensor \si _{1}) = \bigoplus _{j=1} ^p G_j H^0(X,(nL+H+K_X) \tensor \si_0).
\]
\end{thm}

\begin{proof}
By taking $G^{p}= ... G^n = 0$, we may assume that $q = n$.  We are going to use Theorem \ref{main-1} with $F=(n+1)L+H$, $E=L$ and $g^i = G^i$.  Let $\psi = \vp +(n+1)\eta$.  Then 
\[
\ii \di \dbar \psi - qB \ii \di \dbar \eta = \ii \di \dbar \vp + (n+1-qB)\ii \di \dbar \eta \ge 0.
\]
Suppose $f \in H^0(X,\co _X ((n+1)L+H+K_X) \tensor \si_{1})$.  Then 
\[
\int _X \frac{|f|^2e^{-\psi}}{(|g|^2e^{-\eta})^{n+1}}=\int _X \frac{|f|^2e^{-\vp}}{|G|^{2(n+1)}} < +\infty.
\]
By Theorem \ref{main-1} there exist $h_1,...,h_p \in H^0(X,\co _X (nL+H+K_X))$ such that $\sum h_iG^i = f$.  Moreover, 
\begin{eqnarray*}
\int _X \frac{|h|^2e^{-\vp}}{|G|^{2n}R(1-\log |G|^2 +\eta)}&=& \int _X \frac{|h|^2e^{-(\psi-\eta)}}{(|g|^2e^{-\eta})^{2n}R(\xi)} 
\lesssim \int _X \frac{|f|^2e^{-\psi}}{(|g|^2e^{-\eta})^{n+1}}< +\infty.
\end{eqnarray*}
Thus $h_i \in \si _0$ locally.  The proof is complete.
\end{proof}

\section{Hilbert Space Theory}\label{background-section}

\noi {\bf Summation convention.}
We use the the complex version of Einstein's convention, where one sums over (i) repeated indices, one upper and one lower, and (ii) an index and its complex conjugate, provided they are both either upper or lower indices.  In addition, we introduce into our order of operations the rules
\[
|a_i b^i|^2 = |a_1b^1+...|^2, \quad \text{while} \quad 
|a_i|^2|b^i|^2 = |a_1|^2|b^1|^2 + ... \ .
\]

\medskip

\noi {\bf The functional analysis.}
Let $\sh _0$, $\sh _1$, $\sh _2$ and $\sF _1$ be Hilbert spaces with inner products $(\ ,\ )_0$, $(\ ,\ )_1$, $(\ ,\ )_2$ and $(\ ,\ )_*$ respectively.  Suppose we have a bounded linear operator $T_2 : \sh _0 \to \sh _2$ and closed, densely defined operators $T_1 :\sh _0 \to \sh _1$ and $S_1 :\sh _1 \to \sF _1$ satisfying 
$$
S_1T_1 = 0.
$$
Let $\sk = \text{Kernel}(T_1)$.  We consider the following problem.

\begin{pr}\label{u-pr}
Given $\eta \in \sh _2$, is there an element $\xi \in \sk$ such that $T_2 \xi = \eta$?  If so, what can we say about $|\xi| _0$?
\end{pr}

Problem \ref{u-pr} was solved by Skoda in \cite{s-72}.  The difference between Skoda's solution and the one we present here is that Skoda identified the subspace of all $\eta$ for which the problem can be solved, whereas we aim to solve the problems one $\eta$ at a time.  This is a difference in presentation only; the two approaches are equivalent.

\begin{prop}\label{u-soln}
Let $\eta \in \sh _2$.  Suppose there exists a constant $C>0$ such that for all $u \in T_2(\sk)$ and all $\beta \in {\rm Domain}(T_1^*)$, 
\begin{equation}\label{ineq}
|(\eta , u)_2|^2 \le C\left ( |T_2 ^* u +T_1^* \beta |^2_0+|S_1\beta|_*^2\right ).
\end{equation}
Then there exists $\xi \in \sk$ such that $T_2 \xi = \eta $ and $|\xi |_0^2 \le C$.
\end{prop}

\begin{proof}
In \eqref{ineq} we may restrict our attention to $\beta \in {\rm Domain}(T_1^*) \cap \text{Kernel}(S_1)$.  Note that (i) since $S_1T_1 =0$, the image of $T_1 ^*$ agrees with the image of the restriction of $T_1 ^*$ to $\text{Kernel} (S_1)$, and (ii) the image of $T_1 ^*$ is dense in $\sk ^{\perp}$.  Thus the estimate \eqref{ineq} may be rewritten
\begin{equation}\label{q-ineq}
|(\eta ,u)|^2 \le C |[T_2^*u]|^2,
\end{equation}
where we denote by $[\ ]$ the projection to the quotient space $\sh _0 / \sk ^{\perp}$ and the norm on the right hand side is the norm induced on $\sh _0/\sk ^{\perp}$ in the usual way.  As is well known, with this norm $\sh _0/\sk ^{\perp}$ is isomorphic to the closed subspace $\sk$.  (The isomorphism sends any, and thus every, member of $u + \sk ^{\perp}$ to its orthogonal projection onto $\sk $.)  We define a linear functional $\ell : [\text{Image}(T_2^*)]\to \C$ by 
\[
\ell ([T_2 ^*u]) = ( \eta , u)_2.
\]
Then by \eqref{q-ineq} $\ell$ is continuous with norm $\le \sqrt{C}$.  By extending $\ell$ constant in the directions parallel to $[\text{Image} (T_2^*)]^{\perp}$ in $\sh _0 / \sk ^{\perp}$, we may assume that $\ell$ is defined on all of $\sh _0 / \sk ^{\perp}$ with norm still bounded by $\sqrt{C}$.
The Riesz Representation Theorem then tells us that $\ell$ is represented by inner product with respect to some element $\xi$ of $\sh _0 / \sk ^{\perp}$ which we can identify with $\sk$ at this point.  Evidently we have $|\xi|^2_0 \le C$ and 
\begin{equation*}
(T_2 \xi,u)_2 =  (\xi,T_2 ^*u + \sk ^{\perp}) = \ell ([T_2 ^*u]) = (\eta ,u)_2.
\end{equation*}
The proof is complete.
\end{proof}

\medskip

\noi {\bf Hilbert spaces of sections.}
Let $Y$ be a K\"ahler manifold of complex dimension $n$ and $H \to Y$ a holomorphic line bundle equipped with a singular Hermitian metric $e^{-\vp}$.  Given a smooth section $f$ of $H + K_Y\to Y$, we can define its $L^2$-norm
\[
||f||^2 _{\vp} := \int _Y |f|^2 e^{-\vp}.
\]
This norm does not depend on the K\" ahler metric for $Y$.  Indeed, we think of a section of $H+K_Y$ as an $H$-valued $(n,0)$-form.  Then the functions $|f|^2e^{-\vp}$ transforms like the local representatives of a measure on $Y$, and may thus be integrated.

We define 
\[
L^2 (Y,H+K_Y, e^{-\vp})
\]
to be the Hilbert space completion of the space of smooth sections $f$ of $H+K_Y \to Y$ such that $||f||^2 _{\vp} < +\infty$.

More generally, we have Hilbert spaces of $(0,q)$-forms with values in $H+K_Y$.  Given such a $(0,q)$-form $\beta$, defined locally by 
\[
\beta = \beta _{\bar J} d\bar z ^J,
\]
where $J=(j_1,...,j_q) \in \{ 1,...,n\} ^q$ is a multiindex and $d\bar z ^J = d\bar z^{j_1} \wedge ... \wedge d\bar z ^{j_q}$ and the $\beta _{\bar J}$ are skew-symmetric in $J$, we set 
\[
\beta ^J = g^{j_1 \bar k_1} ... g^{j_q \bar k_q} \beta _{\bar K} \quad \text{and} \quad |\beta |^2 = \beta ^J \overline{\beta _{\bar J}},
\]
where $g^{j \bar k}$ is the inverse matrix of the matrix $g_{j \bar k}$ of the K\"ahler metric $g$ of $Y$.  It follows that the functions 
\[
|\beta |^2 e^{-\vp}
\]
transform like the local representatives of a measure, and may thus be integrated.  We then define
\[
||\beta||^2 _{\vp} = \int _Y |\beta|^2 e^{-\vp}.
\]
We now define the Hilbert space 
\[
L^2 _{0,q} (Y,H+K_Y, e^{-\vp})
\]
to be the Hilbert space closure of the space of smooth $(0,q)$-forms $\beta$ with values in $H+K_Y$ such that $||\beta ||^2 _{\vp} < +\infty$.  Of course, these norms depend on the K\" ahler metric $g$ as soon as $q \ge 1$.

We are only going to be (explicitly) interested in the cases $q=0$ and $q=1$, although $q=2$ will enter in an auxiliary way.

\medskip

\noi {\bf Choices.}
In employing Proposition \ref{u-soln}, we shall consider the following spaces.
\begin{eqnarray*}
\sh _0 &:=& ( L^2 (\Omega , K_X + F - E, e^{-\vp_1}))^p\\
\sh _1 &:= & (L^2_{(0,1)}(\Omega ,K_X+F-E,e^{-\vp_1}))^p\\
\sh _2 &:=& L^2(\Omega,K_X+F,e^{-\vp _2})\\
\sF _1 &:= & (L^2_{(0,2)}(\Omega ,K_X+F-E,e^{-\vp_1}))^p
\end{eqnarray*}
Next we define our operators $T_1$ and $T_2$.  Let 
\[
T : L^2 (\Omega , K_X+F-E, e^{-\vp _1} ) \to L^2 _{0,1} (\Omega , K_X+F-E, e^{-\vp _1})
\]
be the densely defined operator whose action on smooth forms with compact support is 
\[
T u = \dbar u.
\]
As usual, the domain of $T$ consists of those $u \in L^2(\Omega,K_X+F-E,e^{-\vp_1})$ such that $\dbar u$, defined in the sense of currents, is represented by an of $L^2_{(0,1)}(\Omega,K_X+F-E,e^{-\vp_1}) $-form with values in $F+K_X$.  We let 
\[
T_1 : \sh _0 \to \sh _1
\]
be defined by 
\[
T_1 (\beta _1,...,\beta _p) = (T\beta _1,...,T\beta _p).
\]
We remind the reader that $T$ has formal adjoint $T^*=T^*_{\vp _1}$ given by the formula
\[
T^*\beta = - e^{\vp _1} \di _{\nu}(e^{-\vp _1} \beta ^{\nu}).
\] 
It follows that
\[
T_1^*(\beta _1,...,\beta _p) = \left ( -e^{\vp _1} \di _{\nu} \left ( e^{- \vp _1}\beta _{1}^{\nu} \right ) , ... ,  -e^{\vp _1}{\di}_{\nu} \left ( e^{- \vp _1}\beta _{p}^{\nu} \right ) \right ).
\]

We will also use the densely defined operators 
\[
S : L^2 _{0,1}(\Omega , K_X+F-E, e^{-\vp _1} ) \to L^2 _{0,2} (\Omega , K_X+F-E, e^{-\vp _1})
\]
defined by $\dbar$ on smooth forms, and the associated operator 
\[
S_1 (\beta _1,..., \beta _p) = (S\beta _1,..., S\beta _p).
\]
However, we will not need the formal adjoint of $S$.

Next we let
\[
T_2 : \sh _0 \to \sh _2
\]
be defined by 
\[
T_2 (h_1,...,h_p) = h_i g^i.
\]
We have
\begin{eqnarray*}
(T_2 ^*u , h)_0 = (u , T_2 h)_2 = \int _{\Omega} u \overline{h_i g^i} e^{-\vp _2} = \int _{\Omega} e^{-(\vp _2- \vp _1)}\bar g ^i u \bar h_j e^{-\vp _1},
\end{eqnarray*}
And thus
\begin{eqnarray*}
T_2 ^* u &=& \left ( e^{-(\vp _2- \vp _1)} \bar g ^1 u , ..., e^{-(\vp _2- \vp _1)} \bar g ^p u \right ).
\end{eqnarray*}

\begin{rmk}
Let us comment on the meaning of this {\it a priori} local formula.  The sections $g ^j$, $1 \le j \le p$, are sections of $E$, and $e^{-(\vp _2 - \vp _1)}$ is a metric for $F - (F-E) = E$.  Thus for each $j$, $|g^j|^2e^{-(\vp _2-\vp _1)}$ is a globally defined function, and so the expressions 
\[
e^{-(\vp _2- \vp _1)} \bar g ^j =e^{-(\vp _2- \vp _1)} |g ^j|^2/g^j
\]  transform like sections of $-E$.  Since $u$ takes values in $K_X +F$, the expressions 
\[
 e^{-(\vp _2- \vp _1)} \bar g ^ju
\]
transform like sections of $K_X + F - E$, which is what we expect.
\end{rmk}

\section{Classical $L^2$ identities and estimates}

In this section we collect some known $L^2$ identities.  

\medskip

\noi {\bf The Bochner-Kodaira Identity.}
Let $\Omega$ be a domain in a complex manifold with smooth, $\R$-codimension-1 boundary $\di \Omega$.  Fix a proper smooth function $\rho$ on a neighborhood of $\Omega$ such that 
\[
\Omega = \{ \rho < 0\} \quad  \text{and} \quad  |\di \rho |\equiv1 \ \text{on } \di \Omega.
\]
Let $H \to \Omega$ be a holomorphic line bundle with singular Hermitian metric $e^{-\vp}$.

The following identity is a basic fact known as the 

\medskip

\noi {\sc Bochner-Kodaira Identity:}  

\noi For any smooth $(0,1)$-form $\beta$ with values in $K_X +H$ and lying in the domain of $\dbar ^*$, 
\begin{eqnarray*}
&& \int _{\Omega} |- e^{\vp} \di_{\nu} (\beta ^{\nu} e^{-\vp})|^2  e^{-\vp} + \int _{\Omega} |\dbar \beta |^2 e^{-\vp} = \int _{\Omega} \beta ^{\nu} \overline{\beta ^{\mu}} ( \di _{\nu} \di _{\bar \mu} \vp ) e^{-\vp}\\
&& \qquad \qquad  +\int _{\Omega} |\bar \nabla \beta|^2 e^{-\vp} + \int _{\di \Omega}   \beta ^{\nu} \overline{\beta ^{\mu}} ( \di _{\nu} \di _{\bar \mu} \rho ) e^{-\vp}.
\end{eqnarray*}

\begin{rmk}
The formal case, in which the boundary term disappears is due to Kodaira, and is a complex version of earlier work of Bochner.  With the boundary term included, the identity above is due to C. B. Morrey.  (For higher degree forms, it is due to Kohn.)
\end{rmk}

\medskip 

\noi {\bf Skoda's Identity.}
For $u \in T_2 (\text{Kernel }T_1)$ and $\beta = (\beta _1,...,\beta _p) \in \text{Domain}({T_1 ^*})$ we have 
\begin{eqnarray*}
&& (T_2 ^* u ,T_1 ^*\beta )_0 \\
&=& (T_1 ( T_2 ^*) u ,\beta) _1 \\
&=& \int _{\Omega} u \left \{ \overline{ \beta _{k} ^{\nu} \di _{\nu} \left ( g^k e^{-(\vp _2-\vp _1)} \right ) } \right \} e^{-\vp _1}.
\end{eqnarray*}
It follows that if $\beta$ is also in $\text{Domain}(S_1)$ then
\begin{eqnarray*}
&& ||T_1^*\beta + T_2 ^* u||_0^2 + ||S _1 \beta ||^2 _* \\
&=& ||T_1^* \beta||_0^2 + ||S_1 \beta ||^2 _* +||T_2 ^* u||_0^2 + 2 \re (T_2 ^* u ,T_1^* \beta )_0 \\
&=& \sum _{k=1} ^p \left ( ||T^* \beta _k||_{\vp _1}^2 + || S \beta _k||^2 _{\vp _1} \right )+ \int _{\Omega} e^{-2 (\vp _2 - \vp _1)} |g|^2 |u|^2 e^{-\vp _1}\\
&&+ 2\re \int _{\Omega}  u \left \{ \overline{ \beta _{k} ^{\nu} \di _{\nu} \left ( g^k e^{-(\vp _2-\vp _1)} \right ) } \right \} e^{-\vp _1}.\\
\end{eqnarray*}
By applying the Bochner-Kodaira identity, we obtain the identity we have called

\noi {\sc Skoda's Identity:}
\begin{eqnarray}\label{skoda-eq}
&&||T_1^*\beta +T_2^* u||^2 _0 + || S_1 \beta ||^2 _* = \int _{\Omega} e^{-(\vp _2-\vp _1)}|g|^2 |u|^2 e^{-\vp _2}\\
\nonumber && \quad + 2\re \int _{\Omega} u \left \{ \overline{\beta _k ^{\nu} \di _{\nu} \left ( g^k e^{-(\vp _2 - \vp _1)}\right )} \right \} e^{-\vp _1} + \int _{\Omega}  ( \beta _k ^{\nu} \overline{\beta _k ^{\mu}} \di _{\nu} \di _{\bar \mu}  \vp _1 ) e^{-\vp _1} \\
\nonumber && \qquad + ||\bar \nabla \beta ||^2_{*}  + \int _{\di \Omega} (\beta _k ^{\nu} \overline{\beta _k ^{\mu}}\di _{\nu} \di _{\bar \mu}\rho) e^{-\vp _1}.
\end{eqnarray}
Here $||\bar \nabla \beta||^2 _* = ||\bar \nabla \beta _1||^2_{\vp _1} + ... + ||\bar \nabla \beta _p||^2_{\vp _1}$.

\medskip 

\noi {\bf Skoda's inequality.}
To obtain an estimate from Skoda's identity, one makes use of the following inequality of Skoda.

\begin{thm}[Skoda's Inequality] \cite{s-72}\label{g-curv-ineq}
Let $g = (g^1,...,g^p)$ be holomorphic functions on a domain $U \subset \C ^n$, and let $q= \min (n,p-1)$.  Then 
\[
q (\beta _k ^{\nu} \overline{\beta _k ^{\mu}}\di _{\nu} \di _{\bar \mu} \log |g|^{2} ) \ge |g|^{2} \left | \beta ^{\nu} _k \di _{\nu}(g^k |g|^{-2}) \right |^2
\]
\end{thm}

\begin{rmk}
When passing to a global setting, it is helpful to keep in mind that $g^k|g|^{-2}$ transforms like a section of the anti-holomorphic line bundle $-\overline E$, and thus $\di _{\nu}(g^k |g|^{-2}))dz ^{\nu}$ transforms like a $(-\overline{E})$-valued $(1,0)$-form.  In particular, both sides of Skoda's inequality consist of globally defined functions.
\end{rmk} 

\medskip

\noi {\bf Skoda's Basic Estimate.}
From Theorem \ref{g-curv-ineq} and Skoda's Identity \eqref{skoda-eq}, we immediately obtain the following slight extension of a theorem of Skoda.

\begin{thm}[Skoda's Basic Estimate]\label{skoda-est-thm}
Let $X$ be an essentially Stein manifold, $E,F \to X$ holomorphic line bundle with singular metrics $e^{-\eta}$ and $e^{-\psi}$ respectively, $\Omega \subset X$ a pseudoconvex domain, $B : \Omega \to (1,\infty)$ a function, and $g^1,...,g^p \in H^0(X,E)$ holomorphic sections.  Set 
\[
q = \min (n,p-1), \quad \vp_1 = \psi + \mu + q\log (|g|^2e^{-\eta})\quad and \quad \vp_2 = \vp_1 + \log |g|^2.
\]
For any $p$-tuple of $F-E$-valued $(0,1)$-forms $\beta = (\beta _1,...,\beta _p) \in {\rm Domain}(T^*_1) \cap {\rm Domain} (S_1)$ and any $u \in T_2 ( {\rm Kernel}(T_1))$ we have the estimate
\begin{eqnarray}\label{skoda-est}
&& ||T^*_1 \beta+T^*_2u||^2 _0 + ||S_1 \beta ||^2_*\\
\nonumber && \quad \ge  \int _{\Omega} \left ( \frac{B-1}{B}\right )|u|^2 e^{-\vp_2} + \int _{\Omega} \tau \beta _k ^{\nu} \overline{\beta _k ^{\mu}}\left ( \di _{\nu} \di _{\bar \mu} \psi - Bq  \di _{\nu} \di _{\bar \mu}\eta \right ) e^{-\vp_1} \\
\nonumber && \qquad + \int _{\Omega} \beta _k ^{\nu} \overline{\beta _k ^{\mu}}\left (\di _{\nu} \di _{\bar \mu} \mu - q(B-1) \di _{\nu} \di _{\bar \mu}\log (|g|^2e^{-\eta} ) \right )e^{-\vp_1}.
\end{eqnarray}
\end{thm}

\begin{proof}
We are going to use Skoda's Identity \eqref{skoda-eq}. First note that, by the Cauchy-Schwartz Inequality, for any open set $U$ we have 
\begin{eqnarray}
\nonumber
&&2\re \int _{U} u \left \{ \overline{ \beta ^{\nu}_k \di _{\nu} \left ( g^k e^{-(\vp _2- \vp _1)}\right )} \right \} e^{-\vp _1} \\
\nonumber
&=&2\re \int _{U} u |g|^{-1} \left \{ \overline{ |g| \beta ^{\nu}_k \di _{\nu} \left ( g^k |g|^{-2}\right )} \right \} e^{-\vp _1} \\
\nonumber
&\ge&  - \int _{U} \frac{1}{B} \frac{|u|^2}{|g|^2} e^{-\vp _1} -  \int _{U} B |g|^{2} \left | \beta _k \di _{\nu} (g^k|g|^{-2})  \right | ^2 e^{-\vp_1} \\
&\ge & - \int _{U}\frac{1}{B} |u|^2 e^{-\vp _2} -  \int _{U}B q (\beta _k ^{\nu} \overline{\beta _k ^{\mu}}\di _{\nu} \di _{\bar \mu} \log |g|^2) e^{-\vp _1}, \label{step}
\end{eqnarray}
where the second inequality follows from Skoda's Inequality (Theorem \ref{g-curv-ineq}).  Since the integrands are globally defined, we may replace $U$ by $\Omega$.  Substituting $\vp _1 = \psi + \mu-q\log |g|^2$, combining the inequality \eqref{step} with Skoda's Identity \eqref{skoda-eq} and dropping the positive terms
\[
|| \bar \nabla \beta||^2 _* \quad \text{and} \quad  \int _{\di \Omega} \beta _k ^{\nu} \overline{\beta _k ^{\mu}}\di _{\nu} \di _{\bar \mu}\rho e^{-\vp _1}
\]
finishes the proof.
\end{proof}

\section{Twisted versions of $L^2$ identities and estimates}

\noi {\bf The twisted Bochner-Kodaira Identity.}
Let $\vp$ be the weight function in the Bochner-Kodaira-H\" ormander identity.  Suppose given a second weight function $\kappa$, and set $\tau = e^{\kappa - \vp}$.  Then 
\[
\ii \di \dbar \vp = \ii \di \dbar \kappa - \frac{\ii \di \dbar \tau }{\tau} + \frac{\ii \di \tau \wedge \dbar \tau }{\tau ^2},
\]
and from the Bochner-Kodaira identity we have
\begin{eqnarray*}
&&\int _{\Omega} |- e^{\kappa} \di _{\nu} ( \beta ^{\nu} e^{-\kappa}) - \tau ^{-1} \beta ^{\nu} \di _{\nu} \tau|^2\tau e^{-\kappa} + \int _{\Omega}\tau |\dbar \beta |^2 e^{-\kappa} \\
&=& \int _{\Omega} \beta ^{\nu} \overline{\beta ^{\mu}} \left ( \tau \di _{\nu} \di _{\bar \mu} \kappa - \di _{\nu} \di _{\bar \mu} \tau +\tau^{-1} (\di _{\nu} \tau)(\di _{\bar \mu} \tau) \right )  e^{-\kappa}\\ 
&& \quad + \int _{\Omega} \tau |\bar \nabla \beta|^2 e^{-\kappa} + \int _{\di \Omega} \tau \beta ^{\nu} \overline{\beta ^{\mu}}\di _{\nu} \di _{\bar \mu}\rho e^{-\kappa}.
\end{eqnarray*} 
Expanding the first term, we obtain the so-called 

\noi {\sc Twisted Bochner-Kodaira Identity}:
\begin{eqnarray*}
||\sqrt{\tau} T^*_{\kappa} \beta||^2_{\kappa} + ||\sqrt{\tau} S \beta||^2_{\kappa} &=& \int _{\Omega} \left ( \tau \di _{\nu} \di _{\bar \mu} \kappa - \di _{\nu}\di _{\bar \mu} \tau \right ) \beta ^{\nu}\overline{\beta ^{\mu}} e^{-\kappa} + 2 \re \int _{\Omega} \beta ^{\nu} ( \di _{\nu}\tau)  \overline{T^* _{\kappa} \beta} e^{-\kappa}\\
&& + \int _{\Omega} \tau |\bar \nabla \beta|^2 e^{-\kappa} + \int _{\di \Omega} \tau \beta ^{\nu} \overline{\beta ^{\mu}}\di _{\nu} \di _{\bar \mu}\rho e^{-\kappa}.
\end{eqnarray*}

\medskip

\noi {\bf Twisted version of Skoda's Identity.}
We shall now twist the weights $\vp _1$ and $\vp _2$ by the same factor $\tau$.  Let
$$
\tau = e^{\kappa _1 - \vp _1} = e^{\kappa _2 - \vp _2}.
$$
Then
$$
\di \dbar \vp _1 = \di \dbar \kappa _1 - \frac{\di \dbar \tau}{\tau} + \frac{\di \tau \wedge \dbar \tau }{\tau ^2},
$$ 
\begin{eqnarray*}
T^* _{\vp _1} \beta &=& - e^{\vp _1} \di _{\nu} ( e^{-\vp _1}\beta ^{\nu}) = - \frac{e^{\kappa _1}}{\tau} \di _{\nu} ( \tau e^{-\kappa _1} \beta ^{\nu}) \\
&=& -e^{\kappa _1}\di _{\nu} (e^{-\kappa _1} \beta ^{\nu}) - \frac{\beta ^{\nu}\di _{\nu} \tau}{\tau} = T^* _{\kappa _1}\beta - \tau  ^{-1} \beta ^{\nu} \di _{\nu} \tau.
\end{eqnarray*}
The operator $T_2^*$ remains unchanged.  We thus calculate that
\begin{eqnarray*}
&& ||T^*_{1,\vp_1}\beta + T^*_2 u ||^2 _{\vp _1}\\
&=& ||\sqrt{\tau} T^*_{1,\kappa_1}\beta - \sqrt{\tau^{-1}} \beta ^{\nu} \di _{\nu} \tau  + \sqrt{\tau} T^*_2 u ||^2_{\kappa_1}\\
&=& ||\sqrt{\tau} T^*_{1,\kappa_1}\beta + \sqrt{\tau} T^*_2 u ||^2_{\kappa_1} - 2 \re \int _{\Omega} (\beta _k ^{\nu} \di _{\nu} \tau) \overline{ (T^*_{1,\kappa _1} \beta + T^*_2 u)_k} e^{-\kappa_1}\\
&& \quad + \int _{\Omega} \tau ^{-1} \beta _k ^{\nu} \overline{\beta _k ^{\mu}}( \di _{\nu} \tau )(\di _{\bar \mu} \tau) e^{-\kappa _1} ,
\end{eqnarray*}
that
\begin{eqnarray*}
\int _{\Omega} u \left \{\overline{\beta _k ^{\nu} \di _{\nu} \left ( g^ke^{-(\vp _2- \vp _1)} \right )}\right \}e^{-\vp _1} =\int _{\Omega} \tau u \left \{\overline{\beta _k ^{\nu} \di _{\nu} \left ( g^ke^{-(\kappa _2- \kappa _1)} \right )}\right \}e^{-\kappa _1}
\end{eqnarray*}
and that 
$$
e^{-\vp _1} \di \dbar \vp _1 = e^{-\kappa _1} (\tau \di \dbar \kappa_1 - \di \dbar \tau - \tau ^{-1} \di \tau \wedge \dbar \tau).
$$
Substitution of these three calculations into Skoda's Identity \eqref{skoda-eq} yields the

\noi {\sc Twisted version of Skoda's Identity:}
\begin{eqnarray}\label{twisted-skoda-eq}
&& ||\sqrt{\tau} T^*_{1,\kappa_1}\beta + \sqrt{\tau} T^*_2 u ||^2_{\kappa_1}+ ||\sqrt{\tau} S_1 \beta ||^2_*= 2\re \int _{\Omega} ( \beta ^{\nu} _k \di _{\nu} \tau) \overline{ (T^*_{1,\kappa _1} \beta + T^*_2 u)_k}e^{-\kappa _1}\\
\nonumber && + \int _{\Omega} e^{-(\kappa _2 - \kappa_1)}|g|^2 \tau |u|^2 e^{-\kappa _2} + 2\re \int _{\Omega} \tau u \left \{ \overline{ \beta ^{\nu} _k \di _{\nu} \left ( g^k e^{-(\kappa _2 -\kappa _1)} \right )}\right \} e^{-\kappa _1}\\
\nonumber && \quad + \int _{\Omega} \beta _k ^{\nu} \overline{\beta _k ^{\mu}}( \tau \di _{\nu} \di _{\bar \mu} \kappa _1 - \di _{\nu} \di _{\bar \mu}\tau )  e^{-\kappa_1}+ \left | \left |\sqrt{\tau}  \bar \nabla \beta \right |\right |^2_* + \int _{\di \Omega} \tau( \beta _k ^{\nu} \overline{\beta _k ^{\mu}}\di _{\nu} \di _{\bar \mu} \rho ) e^{-\kappa _1} .
\end{eqnarray}

\medskip

\noi {\bf Twisted version of Skoda's Basic Estimate.}
The following Lemma is trivial.
\begin{lem}\label{same-domain}
$T_2 ({\rm Kernel}(T_1)) = (T_2 \circ \sqrt{\tau +A}) ( {\rm Kernel}(T_1\circ \sqrt{\tau +A})).$
\end{lem}

Next we have the following result.

\begin{thm}[Twisted version of Skoda's Basic Estimate]\label{twisted-skoda-est}
Let $X$ be an essentially Stein manifold, $E,F \to X$ holomorphic line bundle with singular metrics $e^{-\eta}$ and $e^{-\psi}$ respectively, $\Omega \subset X$ a smoothly bounded pseudoconvex domain, $\tau , B : \Omega \to (1,\infty)$ and $A : \Omega \to (0,\infty)$ functions, and $g^1,...,g^p \in H^0(X,E)$ holomorphic sections.  Set 
\[
q = \min (n,p-1), \quad \kappa _1 = \psi + \mu + q\log (|g|^2e^{-\eta})\quad and \quad \kappa _2 = \kappa _1 + \log |g|^2.
\]
For any $p$-tuple of $F-E$-valued $(0,1)$-forms $\beta = (\beta _1,...,\beta _p) \in {\rm Domain}(T^*_1) \cap {\rm Domain} (S_1)$ and any $u \in (T_2 \circ \sqrt{\tau +A}) ( {\rm Kernel}(T_1\circ \sqrt{\tau +A}))$ we have the estimate
\begin{eqnarray}\label{twisted-skoda-ineq}
&& ||\sqrt{\tau +A} T^*_1 \beta + \sqrt{\tau +A} T^*_2u||^2 _0 + ||\sqrt{\tau} S_1 \beta ||^2_*\\
\nonumber && \quad \ge  \int _{\Omega}\tau \left ( \frac{B-1}{B}\right )|u|^2 e^{-\kappa _2} + \int _{\Omega} \tau \beta _k ^{\nu} \overline{\beta _k ^{\mu}}\left ( \di _{\nu} \di _{\bar \mu} \psi - Bq  \di _{\nu} \di _{\bar \mu}\eta \right ) e^{-\kappa _1} \\
\nonumber && \qquad + \int _{\Omega} \beta _k ^{\nu} \overline{\beta _k ^{\mu}}\left ( \tau \di _{\nu} \di _{\bar \mu} \mu - \di _{\nu}  \di _{\bar \mu} \tau - \frac{(\di _{\nu} \tau  )(\di _{\bar \mu} \tau)}{A} \right )  e^{-\kappa _1} .\\
\nonumber && \qquad +\int _{\Omega} \beta _k ^{\nu} \overline{\beta _k ^{\mu}}\left (- q\tau (B-1) \di _{\nu} \di _{\bar \mu}\log (|g|^2e^{-\eta} ) \right )e^{-\kappa_1}.
\end{eqnarray}
\end{thm}

\begin{proof}
First, by Lemma \ref{same-domain} we may make use of \eqref{twisted-skoda-eq}.

By the pseudoconvexity of $\Omega$, we may drop the last term on the right hand side of \eqref{twisted-skoda-eq}.  The second last term is clearly non-negative and may thus also be dropped.

Since $\kappa _2 - \kappa _1= \log |g|^2$, we see that 
\begin{eqnarray*}
&&2\re \int _{\Omega} \tau u \left \{ \overline{\beta _k ^{\nu} \di _{\nu} \left ( g^k e^{-(\kappa _2 - \kappa _1)}\right )}\right \}e^{-\kappa _1} \\
&\ge& - \int _{\Omega} \frac{\tau}{B} |u|^2 e^{-\kappa _2} - \int _{\Omega} \tau Bq \beta _k ^{\nu} \overline{\beta _k ^{\mu}}( \di _{\nu} \di _{\bar \mu}\log |g|^2)e^{-\kappa _1}.
\end{eqnarray*}
Moreover, by the Cauchy-Schwarz inequality we have 
\begin{eqnarray*}
&&2\re \int _{\Omega} (\beta _k ^{\nu} \di _{\nu} \tau) \overline{(T^*_1  \beta + T_2^*u)_k}e^{-\kappa _1}\\
&\ge& - \int _{\Omega} \frac{\beta _k ^{\nu} \overline{\beta _k ^{\mu}}(\di _{\nu} \tau )(\di _{\bar \mu} \tau )} {A} e^{-\kappa _1} - \int _{\Omega} A|T^*_1  \beta + T_2^*u|^2e^{-\kappa _1}.
\end{eqnarray*}
Applying these inequalities to \eqref{twisted-skoda-eq} easily yields \eqref{twisted-skoda-ineq}.
\end{proof}

\section{Proof of Theorem \ref{main}}\label{main-proof}

We now plan to use Proposition \ref{u-soln} (which solves the first kind of functional analysis problem) to obtain division theorems of type I.  

\medskip

\noi {\bf An a priori estimate.}
Assume that 
\[
\Omega \subset \{ |g|e^{-\eta} < 1\}.
\]
With $q = \min (p-1,n)$, fix a Skoda Triple $(\phi , F,q)$, and let
\[
\xi:= 1-\log (|g|^2e^{-\eta}).
\]
Following Definition \ref{skoda-triple-defn}, let
\[
\tau := \tau (\xi) = \xi+F(\xi) \quad \text{and} \quad A:=A(\xi)=\left ( \frac{-(F''(\xi)+\phi ''(\xi))}{(1+F'(\xi))^2} \right )^{-1}.
\]
Then $\tau \ge 1$.  Let 
\[
\mu = -\phi(\xi).
\]
Then if $\tau$ is non-constant, we have 
\begin{eqnarray*}
&& \tau \ii \di \dbar \mu - \ii \di \dbar \tau - \frac{\ii}{A} \di \tau \wedge \dbar \tau\\
&=& (\tau \phi '(\xi) + 1+F'(\xi))(- \ii \di \dbar \xi) -\left (F''(\xi)+\phi ''(\xi) + \tfrac{(1+F'(\xi))^2}{A}\right ) \ii \di \xi \wedge \dbar \xi\\
&=& (\tau \phi '(\xi) + 1+F'(\xi))(- \ii \di \dbar \xi).
\end{eqnarray*}
The last equality follows from the definition of $A$.  On the other hand, if $\tau$ is constant, we need not apply the twisted Skoda estimate; we can just appeal to the original Skoda estimate, which would yield (by convention)
\begin{eqnarray*}
\tau \ii \di \dbar \mu - \ii \di \dbar \tau - \frac{\ii}{A} \di \tau \wedge \dbar \tau &=& \tau \phi '(\xi) (-\ii \di \dbar \xi) - \phi '' (\xi) \ii \di \xi \wedge \dbar \xi\\
&\ge & (\tau \phi '(\xi)+1+F'(\xi)) (-\ii \di \dbar \xi)
\end{eqnarray*}
Now, $- \di \dbar \xi = \di \dbar \log (|g|^2e^{-\eta})$.  Thus from \eqref{twisted-skoda-ineq} we obtain the {\it a priori} estimate
\begin{eqnarray}\label{twisted-ap}
&&||\sqrt{\tau+A} T^*_1\beta + \sqrt{\tau+A} T^*_2 u ||^2 _{\kappa _1} + ||\sqrt{\tau} S_1 \beta ||^2 _* \\
\nonumber && \quad \ge \int _{\Omega} \tau \frac{B-1}{B}|u|^2 e^{-\kappa _2} + \int _{\Omega} \tau \beta _k ^{\nu} \overline{\beta _k ^{\nu}} ( \di _{\nu}  \di _{\bar \mu} \psi - Bq \di  _{\nu}  \di _{\bar \mu} \eta ) e^{-\kappa _1} \\
\nonumber && \qquad+ \int _{\Omega} \left (-q\tau(B-1) + (\tau \phi ' (\xi) +1+F'(\xi))  \right )\beta _k ^{\nu} \overline{\beta _k ^{\mu}}( \di _{\nu} \di _{\bar \mu}\log (|g|^2e^{-\eta})) e^{-\kappa _1}.
\end{eqnarray}
From the definition of $B$ (See Definition \ref{skoda-triple-defn}) we have 
\[
q\tau(B-1) = \tau \phi ' (\xi)+1+F'(\xi).
\]
Then 
\[
B = \frac{q\tau + \tau \phi ' (\xi) + 1+F'(\xi))}{q\tau},
\]
and thus 
\[
\tau \frac{(B-1)}{B} = \frac{\tau (\tau \phi ' (\xi) + 1+F'(\xi))}{q\tau + \tau \phi ' (\xi) + 1+F'(\xi)}.
\]
Thus we obtain from \eqref{twisted-ap} the estimate
\begin{eqnarray}\label{skoda-triple-est}
&&||\sqrt{\tau+A} T^*_1\beta + \sqrt{\tau+A} T^*_2 u ||^2 _{\kappa _1} + ||\sqrt{\tau} S_1 \beta ||^2 _* \\
\nonumber && \quad \ge \int _{\Omega} \frac{\tau (B-1)}{B} |u|^2 e^{-\kappa _2} + \int _{\Omega} \tau \beta _k ^{\nu} \overline{\beta _k ^{\nu}} ( \di _{\nu}  \di _{\bar \mu} \psi - Bq \di  _{\nu}  \di _{\bar \mu} \eta ) e^{-\kappa _1} .
\end{eqnarray}

\begin{proof}[Conclusion of the proof of Theorem \ref{main}]
Now suppose that 
\[
\ii \di \dbar\psi \ge Bq \ii \di \dbar \eta.
\]
It follows from \eqref{skoda-triple-est} that 
\begin{eqnarray*}
|(f,u)_{\kappa _1}|^2 &\le& \int _{\Omega} \frac{B}{\tau (B-1)}\frac{|f|^2e^{\mu (\tau)} e^{-\psi}}{|g|^{2q+2}} \times \int _{\Omega} \frac{\tau (B-1)}{B} |u|^2 e^{-\kappa _2}
\end{eqnarray*}
In view of Proposition \ref{u-soln}, we find sections $H_1,...,H_p$ such that 
\[
\dbar ( \sqrt{\tau +A} H_i ) = 0 , \quad (g^iH_i)\sqrt{\tau +A} = f
\]
and 
\[
\int _{\Omega} \frac{|H|^2 e^{\phi (\tau)} e^{-\psi}}{|g|^{2q}} \le 
\int _{\Omega} \frac{B}{\tau (B-1)}\frac{|f|^2e^{\phi (\tau)} e^{-\psi}}{|g|^{2q+2}}.
\]
Letting $h _i = \sqrt{\tau +A}H_i$, we obtain 
\[
\int _{\Omega} \frac{|h|^2 e^{\phi (\tau)} e^{-\psi}}{(\tau +A)|g|^{2q}} \le \int _{\Omega} \frac{B}{\tau (B-1)}\frac{|f|^2e^{\phi (\tau)} e^{-\psi}}{|g|^{2q+2}}.
\]
Since the estimates are uniform, we may let $\Omega \to X-V$.  Since $V$ has measure zero, we may replace $X-V$ by $X$.
This completes the proof of Theorem \ref{main}.
\end{proof}

\begin{ack}
I am indebted to Jeff McNeal, from whom I have learned a lot about the method of twisted estimates and with whom I developed the theory of denominators.  I am grateful to Yum-Tong Siu for proposing that I study Skoda's Theorem and try to say something non-trivial about the case $\alpha =1$, which I hope I have done at least somewhat.  
\end{ack}

\end{document}